\documentclass[11pt]{article}

\usepackage{geometry}
\usepackage{amsthm, amsmath, amsfonts, amssymb, mathtools, mathrsfs}
\usepackage{url, hyperref}
\usepackage{xcolor}
\definecolor{lightbrown}{HTML}{C8A27A}
\usepackage{graphicx}

\newtheorem{theorem}{Theorem}[section]
\newtheorem*{theorem*}{Theorem}
\newtheorem{lemma}[theorem]{Lemma}
\newtheorem{proposition}[theorem]{Proposition}

\theoremstyle{remark}

\makeatletter
\newcommand{\subjclass}[2][2020]{%
	\let\@oldtitle\@title%
	\gdef\@title{\@oldtitle\footnotetext{#1 \emph{Mathematics subject classification}: #2}}%
}
\newcommand{\keywords}[1]{%
	\let\@@oldtitle\@title%
	\gdef\@title{\@@oldtitle\footnotetext{\emph{Key words and phrases}: #1}}%
}
\makeatother
\renewcommand{\ge}{\geqslant}
\renewcommand{\le}{\leqslant}

\newcommand{\bbZ}{\mathbb{Z}}
\newcommand{\bbQ}{\mathbb{Q}}
\newcommand{\bbR}{\mathbb{R}}
\newcommand{\bbC}{\mathbb{C}}

\DeclareMathOperator*{\Res}{Res}
\DeclareMathOperator{\rank}{rank}
\DeclareMathOperator{\lcm}{lcm}
\DeclareMathOperator{\sgn}{sgn}
\DeclareMathOperator{\Span}{Span}

\numberwithin{equation}{section}

\allowdisplaybreaks

\newcommand{\coloneq}{:=}

\title{Many more odd zeta values are irrational}
\author{St\'ephane Fischler, Li Lai, Johannes Sprang, and Wadim Zudilin}
\date{}

\subjclass[2020]{Primary 11J72; Secondary 11M06, 33C20.}
\keywords{Odd zeta values, irrationality, hypergeometric series, Hankel determinants.}

\begin{document}
\maketitle


\begin{abstract}
We prove that for any sufficiently large odd integer $s$, there exist at least
\[
2.78 \, \sqrt{\frac{s \log\log s}{\log s}}
\]
irrational values among $\zeta(3), \zeta(5), \dots, \zeta(s)$, where $\zeta$ is the Riemann zeta function. 
This improves on the best known estimate by a constant factor times $\sqrt{ \log\log s}$. 
The proof combines the elimination technique and an arithmetic lemma of Fischler--Sprang--Zudilin, the hypergeometric constructions of Lai--Yu, and the Hankel determinant method that appeared in the proof that $\zeta(5)$ is irrational posted by Fauzan on September 17, 2026.
\end{abstract}

\section{Introduction}

Understanding the arithmetic nature of the values at integers $s\ge 2$ of the Riemann zeta function
\[
\zeta(s)=\sum_{n=1}^\infty \frac1{n^s}
\]
is a longstanding problem. 
When $s$ is even, the zeta value $\zeta(s)$ is a non-zero rational multiple of $\pi^s$: it is transcendental, hence also irrational.

For odd values of $s\ge 3$, no such formula is known (or even conjectured to exist). 
It was only in 1978 when Ap\'ery came up with the proof \cite{Ape1979} of the irrationality of~$\zeta(3)$. 
The next breakthrough in this direction was taken by Ball and Rivoal \cite{BR2001,Riv2000};
they proved that  
\begin{equation}\label{BR}
\dim_{\bbQ} \Span_{\bbQ} \{ 1,\zeta(3),\zeta(5),\ldots,\zeta(s) \} \ge \frac{1-o(1)}{1+\log 2} \,\log s 
\end{equation}
as $s \to \infty$ with $s$ odd.
In particular,
\begin{equation}\label{eqBR}
N_s \ge \frac{1-o(1)}{1+\log 2} \,\log s,
\end{equation}
where 
\[
N_s \coloneq \# \big\{ i \in \{3,5,7,\dots,s\}: \zeta(i) \notin \mathbb{Q} \big\}
\]
counts the amount of irrational numbers among the first zeta values.
The dimension lower bound \eqref{BR} was significantly improved by the first author \cite{Fis2026}; this is however another story.
In this paper, we focus on the improvements of the asymptotic bound \eqref{eqBR}.
It was refined by the first, third and fourth authors \cite{FSZ2018,FSZ2019} using 
an elimination procedure of irrational zeta values (built upon earlier ideas from \cite{Zud2018} and \cite{Spr2018+}):
\begin{equation*} 
N_s \ge 2^{(1-o(1))(\log s) / (\log \log s)}. 
\end{equation*}
The same approach, combined with a more efficient construction, enabled the second author and Yu to obtain the following improvement \cite{LY2020}:
\begin{equation}\label{eqLY}
N_s \ge  1.19 \,\sqrt{\frac{ s}{ \log s}}
\end{equation}
for all sufficiently large odd integers $s$. 
The numerical constant $1.19$ was later replaced with $1.28$ by the second author \cite{Lai2026a}.

Apart from Ap\'ery's theorem, the results mentioned above are asymptotic, namely valid only for sufficiently large~$s$. 
Adapting the proof of \eqref{eqBR}, Rivoal showed \cite{Riv2002} that, among the nine numbers $\zeta(5)$, $\zeta(7)$, \dots, $\zeta(21)$, at least one is irrational. 
This was refined by the fourth author \cite{Zud2001b}: at least one of  $\zeta(5)$, $\zeta(7)$, $\zeta(9)$, $\zeta(11)$ is irrational. 
But the irrationality of $\zeta(5)$ itself was still very open.

On September 17, 2026, a few weeks after powerful AI tools became available, Fauzan posted online a proof that $\zeta(5)$ is irrational \cite{Fau2026+}. 
The principal idea was a new use of Hankel determinants.

In the present paper, we combine such a Hankel determinant with the proof \cite{LY2020} of \eqref{eqLY} to obtain the following improvement.
 
\begin{theorem}\label{main_thm}
For any $\varepsilon > 0$ and all sufficiently large odd integers $s \ge s_0(\varepsilon)$, we have
\[  
N_s \ge \bigg( \sqrt{\frac{4\zeta(2)\zeta(3)}{\zeta(6)}} -\varepsilon \bigg) \cdot \sqrt{\frac{s \log\log s}{\log s}}.  
\]
\end{theorem}

Note that $\sqrt{4\zeta(2)\zeta(3)/\zeta(6)} = 3\sqrt{70\zeta(3)}/\pi^2 = 2.788258\ldots$\,. 
Up to numerical constants, Theorem \ref{main_thm} improves on \eqref{eqLY} by a factor $\sqrt{\log\log s}$.

\medskip
The overall strategy is as follows. 
Assuming that $N_s$ is smaller than the bound of Theorem \ref{main_thm}, we construct (following \cite{LY2020}), for any positive even integer $n$, a linear map $\mathcal{L}_n$ that sends polynomials $P(t) \in \bbQ[t]$ of bounded degree to rational numbers. 
The argument in \cite{LY2020} is to construct such a polynomial $P_n$ in such a way that $\mathcal{L}_n(P_n)\in\mathbb Z\setminus\{0\}$ for any $n$ sufficiently large, and  $\mathcal{L}_n(P_n)$ tends to 0 as $n\to\infty$; this provides a contradiction.
In the present paper, instead of choosing a polynomial $P_n$ we deal with the Hankel determinant
$\Delta_{n,m} \coloneq \det_{0 \le i,j < m}\big( \mathcal{L}_n\left( t^{i+j} \right) \big)\in\mathbb Q$,
where $m$ is not too large in terms of $n$; the details of this construction are given in~\S \ref{sec2}.

The contradiction comes from combining three arguments:
\begin{itemize}
\item An upper bound on the denominator of $\Delta_{n,m}$, obtained in~\S \ref{sec3} by considering a basis consisting of Pochhammer symbols and using standard arithmetic properties. 
No sophisticated $p$-adic estimates are in use. 
\item An upper bound on the absolute value of $\Delta_{n,m}$, obtained in~\S \ref{sec4} also using Pochhammer symbols and elementary estimations.
\item A non-vanishing result that asserts that $\Delta_{n,m}$ is nonzero for at least one value of $m$ in a given range; this is the subject of~\S \ref{sec5}.
\end{itemize}
These results are finally combined in \S \ref{sec6} to complete the proof.

The principal novelties here lie in the non-vanishing part, and our exposition in~\S \ref{sec5} is executed generally enough to serve for future use in related contexts.
It covers Hankel matrices associated with signed measures on discrete closed subsets of the reals.
At the same time our design of the linear functional $\mathcal{L}_n$ deviates significantly from its prototype from~\cite{Fau2026+}.
While using Hankel determinants in \cite{Fau2026+} gives one access to better arithmetic than by traditional techniques, the improvement in Theorem~\ref{main_thm} over the previous bound \eqref{eqLY} comes entirely from the analytic estimates of the determinants.

\medskip
\noindent\textbf{Declaration of generative AI use.}
The initial proof was completely produced by ChatGPT (GPT-6 Astra, OpenAI).
The authors have verified and simplified the argument, rewritten the exposition, and take full responsibility for the mathematical content of this paper.

\medskip
\noindent\textbf{Declaration from natural intelligence.}
The exposition humanly presented in the present note is entirely based on elementary arguments, so that the reader with good knowledge of basics of mathematics can follow them. We also wish to stress that all the ingredients for the methodology in use here have been created before advancements of AI tools. For example, the potentials of Hankel determinant constructions are discussed in \cite{Zud2017} and of their multivariate generalisations in \cite{Bro2026}; the deep connection and a different perspective of the methodology are offered in the framework of holonomicity criteria in \cite{CDT2024+}. Several other generalisations and ingredients are available in related works of these authors and their colleagues. Theorem~\ref{main_thm} has the potential to be generalised to a bound on dimensions of the $\bbQ$-linear spaces generated by odd zeta values rather than the count of irrational numbers among them by implementing a different strategy from \cite{Fis2026}, at least this is something that AI tools suggest; related proofs are expected to be less elementary and more technical though. We see a lot of room for further developments in Diophantine questions about zeta values and, more generally, values of $L$-functions and of polylogarithms. Understanding the beautiful structures behind those methods and proofs seems to us a more exciting task than targeting particular numbers that still await their irrationality status.

\section{Construction} \label{sec2}

Most of the notation, conventions and notions in this section are borrowed from \cite{LY2020};
the new ingredients required for implementation of the Hankel determinant method are the linear functional \eqref{def_L} and determinant~\eqref{def_Delta} itself.

Let $s>100$ be an odd integer, and take a positive real number $B>1$ satisfying $B^2<s/100$.
This number $B$ will be chosen in terms of $s$ in~\S \ref{sec6}.
We will assume that $B$ is larger than some absolute constant so that the estimates in \eqref{est_Z_B} and \eqref{est_A_2} hold.

As in \cite{LY2020}, define the two sets
\begin{gather*}
\Psi_B \coloneq \{ b \in \bbZ_{>0} : \varphi(b) \le B \},
\\
\mathcal{Z}_B \coloneq \Big\{  \frac{a}{b} : b \in \Psi_B, 1 \le a < b \;\text{and}\; \gcd(a,b)=1  \Big\}\subset \bbQ\cap(0,1),
\end{gather*}
where $\varphi(\,\cdot\,)$ is Euler's totient function,
and the integers
\[
A_1(B) \coloneq \prod_{b \in \Psi_B} b^{\varphi(b)}
\quad\text{and}\quad
A_2(B) \coloneq \prod_{b \in \Psi_B} \prod_{\substack{ p \mid b \\ p \;\text{prime}}} p^{\varphi(b)/(p-1)}.
\]

The following rational function is essentially the one from \cite[Definition 2.3]{LY2020} (in which the parameter $r$ is taken to be~0):
for a positive even integer $n$, define
\begin{equation}\label{def_Rn}
R_n(t) \coloneq A_1(B)^n A_2(B)^n n!^{s - |\mathcal{Z}_B|} \cdot \frac{\prod_{\theta \in \mathcal{Z}_B} (t+\theta)_{n}}{(t)_{n+1}^s} \in \bbQ(t),
\end{equation}
where $(t)_m\coloneq\Gamma(t+m)/\Gamma(t)=\prod_{j=0}^{m-1} (t+j)$ denotes the Pochhammer symbol.

The lemmas below are special cases of results proved in \cite{LY2020},
with $r=0$ and $\mathcal{F}_{B} = \mathcal{Z}_{B}\cup\{1\} $ in the notation there.

\begin{lemma}[{\cite[Proposition 2.2]{LY2020}}]
We have the following properties of the sets $\Psi_B$ and $\mathcal{Z}_B$\textup:
\begin{gather}
|\Psi_B|=\bigg(\frac{\zeta(2) \zeta(3)}{\zeta(6)}+o(1)\bigg) B \quad\text{as}\; B\to\infty, \label{est_Psi_B} \\
\Big\{ \frac{a}{b}:1\le a\le b \Big\} \subset \mathcal{Z}_{B}\cup\{1\} \quad\text{for any}\; b \in \Psi_B,
\label{prop_Psi_B} \\
|\mathcal{Z}_B| \le B^2 -1 \quad \text{for any sufficiently large}\; B.
\label{est_Z_B}
\end{gather}
\end{lemma}

\begin{lemma}[{\cite[Lemma 2.4]{LY2020}}]
For the quantities $A_1(B)$ and $A_2(B)$, we have the estimates 
\begin{gather}
A_1(B)=\exp \bigg(\bigg(\frac{\zeta(2) \zeta(3)}{2\zeta(6)}+o(1)\bigg) B^2 \log B\bigg) \quad\text{as}\; B\to\infty,
\label{est_A_1} \\
A_2(B) \le \exp \big(10 B^2(\log \log B)^2\big) \quad\text{for any sufficiently large}\; B. \label{est_A_2}
\end{gather}
\end{lemma}

Our goal is\,---\,until the end of this paper\,---\,to prove that $N_s\ge |\Psi_B|$.
Suppose, for a contradiction, that $N_s<|\Psi_B|$ holds.
Then, we may choose a subset
\[
I \subset \{3,5,\ldots,s\} \quad\text{with}\quad |I| = |\Psi_B|-1,
\]
containing every odd integer $i \in \{3,5,\ldots,s\}$ for which $\zeta(i) \notin \bbQ$.
In other words, $\zeta(i) \in \bbQ$ for all $i \in \{3,5,\ldots,s\} \setminus I$.
Let $A_0 \in \bbZ_{>0}$ be a common denominator of these rational numbers:
\begin{equation}\label{prop_I}
A_0 \zeta(i) \in \bbZ \quad\text{for all}\; i \in \{3,5,\ldots,s\} \setminus I .
\end{equation}
Furthermore, choose and fix integral `weights' $w_b$, where $b \in \Psi_B$, not all zero, such that
\begin{equation}\label{prop_w_b}
\sum_{b \in \Psi_B} w_b b^i =0 \quad\text{for all}\; i \in I. 
\end{equation}
This is clearly possible as the conditions in \eqref{prop_w_b} correspond to a system of homogeneous linear equations in indeterminates $w_b$ with integer coefficients, in which the number of unknowns exceeds the number of equations.

With the choices above, to a polynomial $P(t)\in \bbQ[t]$ assign the number
\begin{equation}\label{def_L}
\mathcal{L}_n(P) \coloneq A_0 \sum_{b \in \Psi_B} w_b \sum_{j=1}^{b} \sum_{\nu=0}^{\infty} 
R_n(t)\Big(P\Big(t+\frac n2\Big)+P\Big(-t-\frac n2\Big)\Big)\Big|_{t=\nu+j/b}.
\end{equation}
This quantity is well defined when $\deg P \le (n+1)s - n|\mathcal{Z}_B| -2$, which ensures that the rational function under summation is $O(t^{-2})$ as $t\to\infty$, hence the series on the right-hand side of~\eqref{def_L} converges absolutely.
In what follows, we will call polynomials that meet the constraint on their degree $n$-admissible.
Finally, for $m \in \bbZ_{>0}$ with $2m \le (n+1)s - n|\mathcal{Z}_B|$, define the Hankel determinant
\begin{equation}\label{def_Delta}
\Delta_{n,m} \coloneq \det_{0 \le i,j < m}\big( \mathcal{L}_n( t^{i+j} ) \big).
\end{equation}

As highlighted in \S 1, the linear functional \eqref{def_L} shows up in \cite{FSZ2019,LY2020} but in a somewhat different light: one makes a choice of appropriate polynomial $P_n(t) \in \bbQ[t]$ (namely, see \cite[Eq.~(2.4) and the beginning of~\S 5]{LY2020}) to demonstrate that $\mathcal{L}_n(P_n)$ is a $\mathbb Z$-linear combination of $1$ and the integers $A_0 \zeta(i) $ in~\eqref{prop_I}, hence an integer, while at the same time it satisfies $0<|\mathcal{L}_n(P_n)|<1$ when $n$~is large.
In a certain sense, the determinant \eqref{def_Delta} gives us access to many such polynomials at once.

\section{Arithmetic estimates}\label{sec3}

Denote $d_n \coloneq \lcm[1,2,\ldots,n]$.
The following lemma follows from the arguments of \cite[Lemma 2]{FSZ2019}.

\begin{lemma}\label{lem_3.1}
Let $n \in 2\bbZ_{>0}$.
Assume that an $n$-admissible polynomial $U(t) \in \bbQ[t]$ satisfies
\[
\frac{d_n^{\lambda}}{\lambda !}\, U^{(\lambda)}(-k) \in \bbZ \quad\text{for any}\; \lambda \in \bbZ_{\ge 0} \;\text{and}\; k \in \bbZ.
\]
Then
$d_n^s \mathcal{L}_n( U ) \in \bbZ$.
\end{lemma}

\begin{proof}
Write $\widetilde{R}_n(t) = R_n(t)(U(t+n/2)+U(-t-n/2)) \in \bbQ(t)$. 
The admissibility of $U(t)$ implies that $\deg \widetilde{R}_n \le -2$, so that the partial-fraction decomposition of $\widetilde{R}_n(t)$ takes the form
\begin{equation}\label{3.1}
\widetilde{R}_n(t) = \sum_{i=1}^{s}\sum_{k=0}^{n} \frac{a_{i,k}}{(t+k)^i}, \quad\text{where}\; a_{i,k} \in \bbQ.
\end{equation}
Since $s$ is odd, $n$ is even, and $\mathcal{Z}_B$ is stable under the reflection $\theta \mapsto 1-\theta$, we have the symmetry $\widetilde{R}_n(-t-n)=-\widetilde{R}_n(t)$, which in particular means that $a_{i,n-k}=(-1)^{i+1}a_{i,k}$ for all $i,k$.
Therefore, for any $b \in \Psi_B$ and $j=1,2,\ldots,b$, the same argument as in \cite[Lemma 1]{FSZ2019} implies that
\begin{equation*}
\sum_{\nu=0}^{\infty} \widetilde{R}_n\Big( \nu+\frac{j}{b} \Big)
= \rho_{0,j/b} + \sum_{\substack{3 \le i \le s \\ i \;\text{odd}}} \rho_{i}\zeta\Big( i,\frac{j}{b} \Big),
\end{equation*}
where $\zeta(i,j/b)$ is the Hurwitz zeta value, its coefficient
\begin{equation}\label{def_rho_i}
\rho_{i} = \sum_{k=0}^{n} a_{i,k} \quad \text{for}\; 3 \le i \le s,\; i \;\text{odd},
\end{equation}
does not depend on $j$, and
\begin{equation}\label{def_rho_0j}
\rho_{0,j/b} = - \sum_{k=1}^{n}\sum_{\ell=0}^{k-1}\sum_{i=1}^{s} \frac{a_{i,k}}{(\ell+j/b)^i}.
\end{equation}
Using the identity $\sum_{j=1}^{b}\zeta(i,j/b) = b^i\zeta(i)$ and Eq.~\eqref{prop_w_b} we obtain
\begin{align}
\mathcal{L}_n(U)
&= A_0 \sum_{b \in \Psi_B} w_b \sum_{j=1}^{b} \sum_{\nu=0}^{\infty} \widetilde{R}_n\Big( \nu+\frac{j}{b} \Big) \notag\\
&= A_0 \sum_{b \in \Psi_B} w_b \sum_{j=1}^{b} \rho_{0,j/b} +  \sum_{\substack{3 \le i \le s \\ i \;\text{odd},\ i \notin I}} \bigg( \sum_{b \in \Psi_B} w_{b}b^i \bigg) \rho_iA_0\zeta(i).
\label{3.4}
\end{align}

Furthermore, for any $i,k$ in \eqref{3.1}, we have
\begin{equation}\label{3.5}
a_{i,k} = \frac{1}{(s-i)!} \big( (t+k)^s\widetilde{R}_n(t) \big)^{(s-i)} \big|_{t=-k}.
\end{equation}
The rational function $(t+k)^s\widetilde{R}_n(t)$ is equal to the product of an integer and the following `building blocks':
$s$ copies of $n!(t+k)/(t)_{n+1}$; the products
\[
b^{n} \cdot \prod_{\substack{p \mid b \\ p \;\text{prime}}} p^{\lfloor n/(p-1) \rfloor} \cdot \frac{(t+a/b)_{n}}{n!},\quad\text{where}\; b \in \Psi_B,\; 1 \le a < b,\; \gcd(a,b)=1;
\]
and the polynomial $U(t+n/2)+U(-t-n/2)$.
By \cite[Lemma 16]{Zud2004b}, \cite[Lemma 4.2]{Lai2026b}, and the condition on $U(t)$, each building block $F(t)$ listed above has the property
\[
\frac{d_n^{\lambda}}{\lambda!} F^{(\lambda)}(-k) \in \bbZ \quad\text{for any $\lambda \in \bbZ_{\ge 0}$ and $k \in \{0,1,\dots,n\}$}. 
\]
Applying the Leibniz rule to \eqref{3.5} and using the property of the building blocks we obtain
\begin{equation}\label{3.6}
d_n^{s-i} a_{i,k} \in \bbZ \quad\text{for all}\; i \in \{1,2,\ldots,s\} \;\text{and}\; k \in \{0,1,\dots,n\}.
\end{equation}
Substituting \eqref{3.6} into \eqref{def_rho_i} and \eqref{def_rho_0j} (for $j=b$) implies that
\begin{equation}\label{3.7}
d_n^{s-i} \rho_i \in \bbZ \quad\text{for any odd $i$ with $3 \le i \le s$}
\end{equation}
and
\begin{equation}\label{3.8}
d_n^{s} \rho_{0,1} \in \bbZ.
\end{equation}

Now, we claim that
\begin{equation}\label{3.9}
d_n^{s} \rho_{0,j/b} \in \bbZ \quad\text{for any $b \in \Psi_B$ and $j \in \{1,2,\ldots,b-1\}$}.
\end{equation}
Indeed, by \eqref{def_rho_0j}, it suffices to prove that
\begin{equation}\label{3.10}
d_{n}^{s} \sum_{i=1}^{s} \frac{a_{i,k_0}}{(\ell_0+j/b)^{i}} \in \bbZ
\end{equation}
for any pair $(k_0,\ell_0)$ with $0 \le \ell_0 < k_0 \le n$.
By \eqref{prop_Psi_B} we know that $j/b \in \mathcal{Z}_B$. 
Thus, by \eqref{def_Rn}, $t=-k_0+\ell_0+j/b$ is a root of $R_n(t)$. 
It follows, after substituting $t=-k_0+\ell_0+j/b$ into \eqref{3.1}, that
\begin{equation*}
\sum_{i=1}^{s} \frac{d_n^{s} a_{i,k_0}}{(\ell_0+j/b)^{i}}
= - \sum_{\substack{0 \le k \le n \\ k \neq k_0}}\sum_{i=1}^{s} \frac{d_n^s a_{i,k}}{(\ell_0-k_0+k+j/b)^{i}}.
\end{equation*}
From this, using the same argument by contradiction as in the part of the proof of \cite[Lemma 2]{FSZ2019} on page 944, we obtain \eqref{3.10} and hence \eqref{3.9}.

Finally, multiplying \eqref{3.4} by $d_n^s$, and using \eqref{3.7}--\eqref{3.9} and \eqref{prop_I}, we obtain the desired integrality of $d_n^s \mathcal{L}_n(U)$.
\end{proof}

We can now provide an explicit denominator for  the Hankel determinant $ \Delta_{n,m} $.

\begin{proposition}\label{prop_ari}
For any $n \in 2\bbZ_{>0}$ and $m \in \bbZ_{>0}$ with $2m\le (n+1)s - n|\mathcal{Z}_B|$, we have
\[
\frac{d_n^{sm}}{\prod_{i=0}^{m-1}\big( n!^{\lfloor i/n \rfloor} \cdot ( i - n\lfloor i/n \rfloor )! \big)^2} \cdot \Delta_{n,m} \in \bbZ.
\]
\end{proposition}

\begin{proof}
Take
\[
U_i(t) = \frac{(t)_i}{n!^{\lfloor i/n \rfloor} \cdot ( i - n\lfloor i/n \rfloor )!} \quad\text{for}\; i=0,1,\dots,m-1.
\]
The polynomial $U_i(t)$ in the collection can be written as the product of the building blocks,
\[
U_i(t)=\frac{(t)_n}{n!}\cdot\frac{(t+n)_n}{n!}\cdot\frac{(t+2n)_{n}}{n!}\dotsb\frac{(t+(\lfloor i/n\rfloor-1) n)_n}{n!}\cdot\frac{(t+n\lfloor i/n\rfloor)_{i - n\lfloor i/n \rfloor}}{( i - n\lfloor i/n \rfloor )!}.
\]
By \cite[Lemma 15]{Zud2004b}, each such a building block $F(t)$ in the product has the property 
\[
\frac{d_n^{\lambda}}{\lambda !} \, F^{(\lambda)}(-k) \in \bbZ \quad\text{for any $\lambda \in \bbZ_{\ge 0}$ and $k \in \bbZ$},
\]
so that by the Leibniz rule the polynomials $U_i(t)U_j(t)$ share this property as well. 
Thus, Lemma \ref{lem_3.1} implies that
\[
d_n^s \mathcal{L}_n( U_iU_j ) \in \bbZ \quad\text{for all}\; 0 \le i,j <m.
\]

The change-of-basis matrix from $1,t,\dots,t^{m-1}$ to $U_0(t),U_1(t),\ldots,U_{m-1}(t)$ is triangular with diagonal entries $1/\big( n!^{\lfloor i/n \rfloor} \cdot ( i - n\lfloor i/n \rfloor )! \big)$.
Therefore,
\[
\det_{0 \le i,j < m}\big( \mathcal{L}_n( U_iU_j ) \big)
= \frac{\Delta_{n,m}}{\prod_{i=0}^{m-1}\big( n!^{\lfloor i/n \rfloor} \cdot ( i - n\lfloor i/n \rfloor )! \big)^2},
\]
hence
\[
\frac{d_n^{sm}}{\prod_{i=0}^{m-1}\big( n!^{\lfloor i/n \rfloor} \cdot ( i - n\lfloor i/n \rfloor )! \big)^2} \cdot \Delta_{n,m}
= \det_{0 \le i,j < m}\big( d_n^s \mathcal{L}_n( U_iU_j ) \big)\in\bbZ.
\qedhere
\]
\end{proof}

\section{Analytic behaviour}\label{sec4}

The following proposition gives an upper bound for the absolute value of the Hankel determinant $ \Delta_{n,m}$ with the normalisation factor from Proposition~\ref{prop_ari}.

\begin{proposition}\label{prop_ana}
For any $n \in 2\bbZ_{>0}$ and $m \in \bbZ_{>0}$ such that
\begin{equation}\label{ana_hypo}
2m \le (s-B^2-10)n,
\end{equation}
we have
\begin{multline}
\frac{d_n^{sm}}{\prod_{i=0}^{m-1}\big( n!^{\lfloor i/n \rfloor} \cdot ( i - n\lfloor i/n \rfloor )! \big)^2} \cdot |\Delta_{n,m}| \\
\le(C_0m)^m A_1(B)^{mn}A_2(B)^{mn} e^{m^2+(2s+10B^2)mn} \cdot \bigg( \frac{10B^2n}{m} \bigg)^{(s-B^2)mn-2m^2},
\label{ana_conc}
\end{multline}
where the constant $C_0$ is independent of $m$ and~$n$.
\end{proposition}

\begin{proof}
We will use a crude form $d_n \le e^{2n}$ of Chebyshev's bound for Chebyshev's $\psi$ function and the elementary inequality
\begin{equation}\label{weak_Stirling}
\log N! = \sum_{j=1}^N \log j \ge \int_1^N\log x\,\mathrm{d}x \ge N\log N-N \quad\text{for all}\; N \in \bbZ_{>0},
\end{equation}
which remains true for $N=0$ with the convention $0\log 0 = 0$.
For any $i=0,1,\dots,m-1$, from the last inequality and $x\log x \ge -e^{-1} \ge -1$ for $x \in [0,1]$ we obtain 
\begin{align*}
\log\big( n!^{\lfloor i/n \rfloor} \cdot ( i - n\lfloor i/n \rfloor )! \big)
&= \lfloor i/n \rfloor \log n! + \log ( n\{ i/n \} )! \\
&\ge \lfloor i/n \rfloor ( n\log n - n ) +  n\{ i/n \} \log(  n\{ i/n \} ) -  n\{ i/n \} \\
&= i \log n - i + n\{ i/n \} \log \{ i/n \} \\
&\ge i\log n - i -n.
\end{align*}
Therefore,
\begin{align}
\frac{d_n^{sm}}{\prod_{i=0}^{m-1}\big( n!^{\lfloor i/n \rfloor} \cdot ( i - n\lfloor i/n \rfloor )! \big)^2}
&\le e^{2smn}\cdot\exp\bigg( 2 \sum_{i=0}^{m-1} ( -i\log n + i + n ) \bigg) \notag\\
&\le e^{2smn -m(m-1)\log n + m^2 + 2mn} = \frac{e^{m^2+2(s+1)mn}}{n^{m(m-1)}}. \label{4.3}
\end{align}

Now, for any $i=0,1,\dots,m-1$, consider the polynomial
\begin{equation}\label{def_Vi}
V_i(t) = t^{i - ( 2|\mathcal{Z}_B| + 2 )\lfloor i/( 2|\mathcal{Z}_B| + 2 ) \rfloor} \cdot \prod_{\substack{\theta \in \mathcal{Z}_B \cup \{1\} \\ 0 \le \nu \le \lfloor i/\left( 2|\mathcal{Z}_B| + 2 \right) \rfloor - 1}} \Big( t^2 - \Big( \nu + \theta + \frac{n}{2} \Big)^2 \Big),
\end{equation}
where the empty product is understood as $1$. 
Note that $V_i$ is a monic polynomial of degree~$i$. 
Therefore, the change-of-basis matrix from $1,t,\ldots,t^{m-1}$ to $V_0(t),V_1(t),\ldots,V_{m-1}(t)$ is unitriangular; in particular,
\begin{equation}\label{4.4}
\Delta_{n,m} = \det_{0 \le i,j < m}\big( \mathcal{L}_n(V_iV_j) \big).
\end{equation}
Furthermore, $\left|\mathcal{L}_n(V_i^2)\right|\le S_i$ for $i=0,1,\dots,m-1$, where
\begin{equation}\label{def_Si}
S_i \coloneq A_0 \sum_{b \in \Psi_B} |w_b| \sum_{j=1}^{b} \sum_{\nu=0}^{\infty} \Big| R_n\Big( \nu+\frac{j}{b} \Big) \Big| \Big( \Big| V_i\Big(\nu+\frac{j}{b}+\frac{n}{2}\Big) \Big|^2 + \Big| V_i\Big(-\nu-\frac{j}{b}-\frac{n}{2}\Big) \Big|^2 \Big);
\end{equation}
the series on the right-hand side converges absolutely because \eqref{ana_hypo} and \eqref{est_Z_B} imply
\[
\deg R_n + 2\deg V_i \le n|\mathcal{Z}_B| - (n+1)s + 2(m-1) < -2.
\]
By \eqref{def_L} and the Cauchy--Schwarz inequality we have
$| \mathcal{L}_n(V_iV_j)  | \le \sqrt{S_iS_j}$,
in turn implying that
\begin{equation}\label{4.6}
| \Delta_{n,m} |
\le \sum_{\sigma \in \mathfrak{S}_m} \prod_{i=0}^{m-1} | \mathcal{L}_n( V_iV_{\sigma(i)} ) |
\le m! \prod_{i=0}^{m-1} S_i
\le m^m \prod_{i=0}^{m-1} S_i.
\end{equation}

To estimate the individual terms in $S_i$, consider $b,j,\nu$ in the summation range of \eqref{def_Si} and $i=0,1,\dots,m-1$. 
By \eqref{prop_Psi_B} and \eqref{def_Vi}, we have
\begin{alignat}{3}
V_i\Big(\nu+\frac{j}{b}+\frac{n}{2}\Big)
&=  V_i\left(-\nu-\frac{j}{b}-\frac{n}{2}\right)
&&= 0 
&\quad&\text{if}\; \nu < \lfloor i/( 2|\mathcal{Z}_B| + 2 ) \rfloor,
\label{4.9}\\
\Big| V_i\Big(\nu+\frac{j}{b}+\frac{n}{2}\Big) \Big|
&= \Big| V_i\Big(-\nu-\frac{j}{b}-\frac{n}{2}\Big) \Big|
&&\le \Big( \nu+\frac{j}{b}+\frac{n}{2}\Big)^i  \notag\\ 
&&&\le \Big( \nu+\frac{j}{b}+n\Big)^i 
&\quad&\text{if}\; \nu \ge \lfloor i/( 2|\mathcal{Z}_B| + 2 ) \rfloor.
\label{4.10}
\end{alignat}
For any real number $t$ of the form $t = \nu+j/b$, we have $t \ge 1/(B_0-1)$, where $B_0=1+\max_{b\in\Psi_B}b\in\bbZ_{>1}$.
Therefore, 
\begin{align*}
(t)_{n+1} = n! \cdot t \cdot \prod_{k=1}^{n}\left( 1+\frac{t}{k} \right) &\ge n! \cdot t \cdot \left( 1+\frac{t}{n} \right)^n \\ 
&= n! \cdot \left( 1+\frac{t}{n} \right)^{n+1} \cdot \frac{1}{\frac{1}{t}+\frac{1}{n}} \ge n! \cdot \left( 1+\frac{t}{n} \right)^{n+1} \cdot \frac{1}{B_0}.
\end{align*}
Obviously, $(t+\theta)_{n} \le (t+n)^{n}$ for all $\theta \in \mathcal{Z}_B \subset (0,1)$.
Substituting these two estimates into \eqref{def_Rn} with $t=\nu+j/b$ and using \eqref{weak_Stirling} we obtain
\begin{align}
\Big| R_n\Big( \nu+\frac{j}{b} \Big) \Big|
&\le B_0^s \cdot A_1(B)^nA_2(B)^n \cdot \frac{n^{(n+1)s}}{n!^{|\mathcal{Z}_B|}} \cdot \frac{1}{\left(\nu+j/b+n\right)^{(n+1)s-n|\mathcal{Z}_B|}}
\notag\\
&\le B_0^s \cdot A_1(B)^nA_2(B)^n e^{|\mathcal{Z}_B|n} \cdot \frac{n^{(n+1)s-n|\mathcal{Z}_B|}}{\left(\nu+j/b+n\right)^{(n+1)s-n|\mathcal{Z}_B|}}.
\label{4.13}
\end{align}

Now, substituting \eqref{4.9}, \eqref{4.10} and \eqref{4.13} into \eqref{def_Si} we deduce that
\begin{align}
S_i
&\le 2A_0 \cdot B_0^s \cdot A_1(B)^n A_2(B)^ne^{|\mathcal{Z}_B|n} \cdot n^{2i}
\notag\\ &\qquad\times
\sum_{b \in \Psi_B} |w_b| \sum_{j=1}^{b} \sum_{\nu=\lfloor i/( 2|\mathcal{Z}_B| + 2 ) \rfloor}^{\infty} \frac{n^{(n+1)s-n|\mathcal{Z}_B|-2i}}{(\nu+j/b+n)^{(n+1)s-n|\mathcal{Z}_B|-2i}}. \label{4.16}
\end{align}
For brevity, write $\nu_i = \lfloor i/\left( 2|\mathcal{Z}_B| + 2 \right) \rfloor$ and $e_i = (n+1)s-n|\mathcal{Z}_B|-2i$ for $i=0,1,\ldots,m-1$. Then
\begin{align}
\sum_{\nu=\nu_i}^{\infty} \frac{n^{e_i}}{(\nu+j/b+n)^{e_i}}
&\le \frac{n^{e_i}}{( \nu_i + n )^{e_i}} + \int_{\nu_i}^{\infty} \frac{n^{e_i}}{(x+n)^{e_i}} \,\mathrm{d}x \notag\\
&=\frac{1}{( 1 + \nu_i/n )^{e_i-1}}\bigg( \frac{1}{1+ \nu_i/n} + \frac{n}{e_i-1} \bigg)
\le \frac{2}{( 1 +  \nu_i/n )^{e_i-1}}, \label{4.17}
\end{align}
where the last step used $e_i-1 > 10n$ by \eqref{ana_hypo}.
With the bound \eqref{est_Z_B} on $|\mathcal{Z}_B|$, we have
\[
1+\frac{\nu_i}{n} > 1+\frac{1}{n}\bigg( \frac{i}{2B^2} -1 \bigg) > \frac{i+1}{2B^2n}.
\]
Combining this with $e_i-1 > (n+1)s-n|\mathcal{Z}_B| -2m > (s-B^2)n-2m$ and \eqref{4.17} continues the estimate \eqref{4.16} as follows:
\begin{equation}\label{4.19}
S_i \le C_0 \cdot A_1(B)^nA_2(B)^n e^{B^2 n} \cdot n^{2i} \cdot \bigg( \frac{2B^2 n}{i+1} \bigg)^{(s-B^2)n-2m},
\end{equation}
where $C_0 =  4A_0 \cdot B_0^s \sum_{b \in \Psi_B} b|w_b|$ is independent of $m$ and~$n$.
Substituting \eqref{4.19} into \eqref{4.6}, and using \eqref{weak_Stirling}, we obtain
\begin{align}\label{4.20}
| \Delta_{n,m} |
&\le (C_0m)^m  A_1(B)^{mn}A_2(B)^{mn} e^{B^2 mn} \cdot n^{m(m-1)} \cdot \bigg( \frac{(2B^2n)^m}{m!} \bigg)^{(s-B^2)n-2m} \notag\\
&\le (C_0m)^m  A_1(B)^{mn}A_2(B)^{mn} e^{B^2 mn} \cdot n^{m(m-1)} \cdot \bigg( \frac{2eB^2n}{m} \bigg)^{(s-B^2)mn-2m^2}.
\end{align}
Finally, putting \eqref{4.3} and \eqref{4.20} together, the desired estimate \eqref{ana_conc} follows.
\end{proof}

\section{Non-vanishing}\label{sec5}
In this section we prove a non-vanishing result for the Hankel determinants $\Delta_{n,m}$. 
The argument works more generally for a wider class of moment matrices of signed measures supported on discrete closed subsets $\mathcal{X}\subseteq \bbR$. 
We first give a general non-vanishing result for such Hankel matrices of moments and then apply it to the determinants $\Delta_{n,m}$.

\subsection{Hankel determinants of moments of discrete measures}\label{subsec:Hankel-det}
Let $\mathcal{X}\subseteq \bbR$ be a closed discrete subset and $\mu$ a signed measure on $\mathcal{X}$. 
We write $(a_x)_{x\in \mathcal{X}}$ for the sequence of weights of $\mu$:
\[
\mu=\sum_{x\in \mathcal{X}} a_x\delta_x,
\]
where $\delta_x$ denotes the delta distribution at $x$. 
We assume that the sequence satisfies a growth condition of the form
\begin{equation}\label{eq:growth-condition-GD}
\sum_{x\in \mathcal{X}} |a_x|(1+|x|)^{D-1}<+\infty\tag{$\mathrm{G}_D$}
\end{equation}
for some positive integer $D$. 
This growth condition implies that the \emph{Cauchy--Stieltjes transform} of $\mu$ is defined as the meromorphic function
\begin{equation}\label{eq:def-CS}
f(z):=\int_{\mathcal{X}} \frac{\mathrm{d}\mu(t)}{z-t}=\sum_{x\in \mathcal{X}} \frac{a_x}{z-x}.
\end{equation}
Indeed, the growth condition guarantees that the right hand side converges absolutely and uniformly on every compact subset of $\bbC\setminus \mathcal{X}$. 
By definition, $f$ has (at most) simple poles at $x\in \mathcal{X}$ with residues
\[
\Res_{z=x} f(z)=a_x.
\] 
In the following, it is convenient to extend the weights by zero to all real numbers, i.e., we define
\[
a(t):=\begin{cases}
a_t & t\in \mathcal{X},\\
0 & t\in \bbR\setminus\mathcal{X}.
\end{cases}
\]
The following lemma shows that the Cauchy--Stieltjes transform $f$ of $\mu$ satisfies a functional equation with polynomial coefficients if the weight function satisfies an identity of the form $A(t)a(t+1)=B(t)a(t)$ with polynomials $A,B\in \bbC[t]$, compare \cite[Theorem 11]{DM2014}:
\begin{lemma}\label{lem:feq-CS-transform}
Let $\mu$ be a signed measure and $D\ge 1$ be a positive integer such that the weights of $\mu$ satisfy the growth condition \eqref{eq:growth-condition-GD}. 
Suppose there exist polynomials $A(z),B(z)\in \bbC[z]$ of degree $\le D$ such that $A(t)a(t+1)=B(t)a(t)$ for all $t\in \bbR$. Then the Cauchy--Stieltjes transform $f$ of~$\mu$ satisfies the functional equation
\begin{equation}\label{eq:CS-feq}
A(z)f(z+1)-B(z)f(z)=C(z),
\end{equation}
where $C(z)\in \bbC[z]$ is the polynomial of degree $\le \max(\deg(A),\deg(B))-1$ (with the convention $C=0$ if $\deg A=\deg B=0$) given by
\begin{equation}\label{eq:def-C}
C(z)=\sum_{x\in \mathcal{X}}a_x\left(\frac{A(z)-A(x-1)}{z-(x-1)}- \frac{B(z)-B(x)}{z-x}\right).
\end{equation}
\end{lemma}

\begin{proof}
First, note that each of the terms
\[
\frac{A(z)-A(x-1)}{z-(x-1)} \quad\text{and}\quad \frac{B(z)-B(x)}{z-x}
\]
in the definition of $C(z)$ is a polynomial of degree $\le \max(\deg(A),\deg(B))-1$ (with the convention $C=0$ if $\deg A=\deg B=0$). 
The growth condition \eqref{eq:growth-condition-GD} guarantees that the coefficients in the formula defining $C(z)$ converge absolutely. 
Using
\[
\frac{A(z)}{z-(x-1)}=\frac{A(z)-A(x-1)}{z-(x-1)}+\frac{A(x-1)}{z-(x-1)},
\]
and the corresponding identity for $B(z)$, we obtain
\begin{alignat*}{2}
\kern18.5mm
A(z)f(z+1)-B(z)f(z)-C(z)
&=\sum_{x\in \mathcal{X}}\frac{a_x A(x-1)}{z-(x-1)}-\sum_{x\in \mathcal{X}}\frac{a_xB(x)}{z-x}
\\
&= \sum_{t\in \mathcal{X}\cup (\mathcal{X}-1)} \frac{A(t)a(t+1)-B(t)a(t)}{z-t}=0.
&\kern30mm\qedhere
\end{alignat*}
\end{proof}

The following lemma says that any sufficiently good rational approximation of a Cauchy--Stieltjes function $f$ with a functional equation of the form \eqref{eq:CS-feq} satisfies the same functional equation:

\begin{lemma}\label{lem:feq-for-good-approximations}
Let $f$ be the Cauchy--Stieltjes transform of a signed measure $\mu$ on $\mathcal{X}$, satisfying a functional equation of the form
\[
A(z)f(z+1)-B(z)f(z)=C(z),
\]
with $A(z),B(z),C(z)\in \bbC[z]$. 
Write $d:=\max(\deg A(z),\deg B(z))$.
Suppose there exist polynomials $Q(z),P(z)\in \bbC[z]$ with $Q \neq 0$ such that
\begin{equation}\label{eq:approx-f}
Q(z)f(z) - P(z) = o\left(y^{-d-\deg Q(z)}\right) \quad\text{as}\; z=\mathbf{i}y \;\text{or}\; z=1+\mathbf{i}y \;\text{and}\; y \to \infty.
\end{equation}
Then the rational function $r(z):=P(z)/Q(z)$ satisfies the same functional equation as $f(z)$:
\[
A(z)r(z+1)-B(z)r(z)=C(z).
\]
\end{lemma}

\begin{proof}
Consider the polynomial
\begin{equation}\label{def_W}
W(z) = A(z)P(z+1)Q(z) - B(z)P(z)Q(z+1) - C(z)Q(z)Q(z+1) \in \bbC[z].
\end{equation}
Using the functional equation of $f(z)$ we can rewrite this as
\begin{align*}
W(z) &= A(z)Q(z)\big( P(z+1) - Q(z+1)f(z+1) \big) - B(z)Q(z+1) \big( P(z) - Q(z)f(z) \big).
\end{align*}
The assumption \eqref{eq:approx-f} gives
\[
W(z)= o(1)\quad\text{as $z=\mathbf{i}y \to \mathbf{i}\infty$}. 
\]
Thus, the polynomial $W(z)$ vanishes identically.
Dividing \eqref{def_W} by $Q(z)Q(z+1)$ then gives
\[
A(z) r(z+1)- B(z) r(z)= C(z).
\qedhere
\]
\end{proof}

The following is a classical result about Hankel matrices in the framework of Kronecker's theory of recurrent sequences.

\begin{lemma}\label{lem_window}
Let $M,N \in \bbZ_{>0}$ with $ M \le N$, and let $\mu_0,\ldots,\mu_{2N-2}$ belong to a field. 
If
\[
\rank (\mu_{i+j})_{\substack{0\le i < N \\ 0 \le j <M}}=M,
\]
then $\det_{0 \le i,j<m} (\mu_{i+j})\ne 0$ for some $M \le m \le N$.
\end{lemma}

\begin{proof}
In the square Hankel matrix $(\mu_{i+j})_{0 \le i,j <N}$, let $m$ be the largest integer such that its first $m$ columns are linearly independent. 
Then $m \ge M$ by the condition. 
If $m=N$, the conclusion is immediate. 
Otherwise, its first $m+1$ columns are dependent, and \cite[Lemma~2.3]{ANL2021} gives $\det_{0 \le i,j < m}(\mu_{i+j}) \ne 0$.
\end{proof}

Suppose that $\mu$ satisfies the growth condition \eqref{eq:growth-condition-GD}, then 
for any $0\le k<D$, we define its $k$-th \emph{moment} by
\[
\mathcal{L}(t^k):=\int_{\mathcal{X}} t^k \mathrm{d}\mu(t)=\sum_{x\in \mathcal{X}} a_x x^k.
\]
In particular, the Hankel determinant
\[
\Delta_m:=\det_{0\le i,j<m}( \mathcal{L}(t^{i+j}))
\]
of the moment sequence is defined for $2m-2<D$. 
The following proposition gives a non-vanishing criterion for these Hankel determinants.

\begin{proposition}\label{prop:nonvanishing-Hankel-det}
Let $\mu$ be a signed measure on $\mathcal{X}$ satisfying the growth condition \eqref{eq:growth-condition-GD} for some positive integer~$D$. 
Assume, in addition, that the extension by zero $a(t)$ of weights satisfies
\[
A(t)a(t+1)=B(t)a(t)
\]
for some polynomials $A(z),B(z)\in \bbC[z]$ with $1 \le d:=\max(\deg A(z),\deg B(z))\le D$ and that the equation
\[
A(z)r(z+1)-B(z)r(z)=C(z)
\]
has no solution in rational functions $r(z)\in \bbC(z)$, where
\[
C(z)=\sum_{x\in \mathcal{X}}a_x\left(\frac{A(z)-A(x-1)}{z-(x-1)}- \frac{B(z)-B(x)}{z-x}\right)
\] 
is defined as in \eqref{eq:def-C}. 
Then for any integer $M\ge 1$ such that $2(M+d-1)-2<D$ there exists a positive integer $m$ such that
\[
M\le m\le M+d-1 \quad \text{and}\quad \Delta_m\neq 0,
\]
where $\Delta_m$ is the $m$-th Hankel determinant of moments of~$\mu$.
\end{proposition}

\begin{proof}
Put $N=M+d-1$. 
Since $2N-2<D$, the Hankel determinants $\Delta_m$ are defined for all $m\le N$. 
Let us write $f(z)$ for the Cauchy--Stieltjes transform of~$\mu$. 
According to Lemma \ref{lem:feq-CS-transform}, $f(z)$ satisfies the functional equation
\[
A(z)f(z+1)-B(z)f(z)=C(z)
\]
with $A(z),B(z),C(z)$ as in the statement of the proposition.
By Lemma \ref{lem_window}, it suffices to prove that
\begin{equation}\label{eq:rk-Hankel-mat}
\rank \big(\mathcal{L}(t^{i+j})\big)_{\substack{0\le i < N \\ 0 \le j <M}}=M.
\end{equation}
We argue by contradiction. 
More precisely, we will see that failure of \eqref{eq:rk-Hankel-mat} implies the existence of a good rational approximation which will automatically satisfy the same functional equation as $f(z)$, contradicting the assumptions of the proposition. 
Let us now start the proof by contradiction. 
Suppose \eqref{eq:rk-Hankel-mat} fails.
Then, there exists a column vector $(c_0,c_1,\dots,c_{M-1})^\top \in \operatorname{Mat}_{M \times 1}(\bbC) \setminus \{\boldsymbol{0}\}$ such that
\begin{equation}\label{555}
\big(\mathcal{L}(t^{i+j})\big)_{\substack{0\le i < N \\ 0 \le j <M}} \cdot (c_0,c_1,\ldots,c_{M-1})^\top = \boldsymbol{0}. 
\end{equation}
Write 
\[
Q(t)=\sum_{j=0}^{M-1}c_jt^j \in \bbC[t] \setminus \{0\}, \quad \deg Q \le M-1.
\] 
Then \eqref{555} implies that
\begin{equation}\label{5.12}
\sum_{x \in \mathcal{X}} a_x Q(x)x^{i} = 0 \quad\text{for all}\; i=0,1,\ldots,N-1.	
\end{equation}
Define
\[
P(z) = \sum_{x \in \mathcal{X}} a_x \frac{Q(z)-Q(x)}{z-x}.
\]
The growth condition \eqref{eq:growth-condition-GD} ensures that $P(z) \in \bbC[z]$ with $\deg P \le M-2$ (if $M=1$, then $P=0$). 
Then
\begin{align*}
Q(z)f(z) - P(z) &= \sum_{x \in \mathcal{X}} a_x Q(x) \frac{1}{z-x} \quad\text{by \eqref{eq:def-CS}} \\
&= \sum_{x \in \mathcal{X}} a_x Q(x) \bigg( \frac{x^N}{z^N(z-x)} + \sum_{i=0}^{N-1}\frac{x^i}{z^{i+1}} \bigg) \\
&= \sum_{x \in \mathcal{X}} a_x Q(x) \frac{x^N}{z^{N}(z-x)} \quad\text{by \eqref{5.12}}\\
&= \frac{1}{z^N}\sum_{x \in \mathcal{X}} a_x Q(x)x^{N-1} \frac{x}{z-x}. 
\end{align*}
The growth condition \eqref{eq:growth-condition-GD} ensures $\sum_{x \in \mathcal{X}} |a_xQ(x)x^{N-1}| < \infty$, since $\deg Q + N-1 \le 2N-2 < D$.
For $z=\mathbf{i}y$ or $z=1+\mathbf{i}y$ with $y \ge 1$, the quotient $x/(z-x)$ has absolute value $\le 2$ and tends to $0$ as $y \to \infty$.
Therefore, the dominated convergence theorem implies
\begin{equation}\label{5.13}
Q(z)f(z) - P(z) = o(y^{-N}) \quad\text{as}\; z=\mathbf{i}y \;\text{or}\; z=1+\mathbf{i}y \;\text{and}\; y \to \infty.
\end{equation}
Since $N=M+d-1\ge \max(\deg A(z),\deg B(z))+\deg Q(z)$, Lemma \ref{lem:feq-for-good-approximations} shows that the rational function $r(z)=P(z)/Q(z)$ satisfies the same functional equation as $f(z)$, i.e.
\begin{equation*}
A(z)\frac{P(z+1)}{Q(z+1)}-B(z)\frac{P(z)}{Q(z)}=C(z),
\end{equation*}
contradicting the assumption in the proposition.
\end{proof}

\subsection{Application}
For convenience, when $n \in 2\bbZ_{>0}$, define the shift of the rational function~\eqref{def_Rn},
\[
\widehat{R}_n(t) \coloneq R_n\Big( t-\frac{n}{2} \Big) \in \bbQ(t).
\]
Since $n$ is even, $s$ is odd, and $\mathcal{Z}_B$ is stable under the reflection $\theta \mapsto 1-\theta$, we can rephrase the (well-poised) property $R_n(-t-n) = -R_n(t)$ as
\begin{equation}\label{oddness_of_shifted_Rn}
\widehat{R}_n( -t ) = - \widehat{R}_n(t).
\end{equation}
Define the set
\[
\mathcal{X}_B \coloneq \bbZ + (\mathcal{Z}_B \cup \{1\})
= \big\{  \nu+\theta : \nu \in \bbZ,\ \theta \in \mathcal{Z}_B \cup \{1\} \big\} \subset \bbQ
\]
and the weights
\[
\widehat{w}_{x} \coloneq \sum_{\substack{b \in \Psi_B \\ bx \in \bbZ}} w_b \in \bbZ
\quad\text{for}\; x \in \mathcal{X}_B,
\]
so that
\begin{equation*} 
\widehat{w}_{x+1} = \widehat{w}_x = \widehat{w}_{-x} \quad\text{for}\; x \in \mathcal{X}_B.
\end{equation*}
We want to apply the theory in \S\ref{subsec:Hankel-det} to prove the non-vanishing of the determinants $\Delta_{n,m}$ introduced in \eqref{def_Delta}. 
For this, we introduce the discrete measure
\[
\mu_n:=A_0\sum_{x \in \mathcal{X}_B}  a(x) \delta_x = A_0\sum_{x \in \bbR}  a(x) \delta_x
\]
assigned to
\[
a(x):=\begin{cases}
\sgn(x)\widehat{w}_x \widehat{R}_n(x) &\text{for}\; x\in \mathcal{X}_B \;\text{and}\; |x|>n/2,\\
0 & \text{otherwise}.
\end{cases}	
\]
With the help of \eqref{def_Rn}, we get
\begin{equation}\label{5.8}
\frac{\widehat{R}_n( t + 1)}{\widehat{R}_n( t )}
=\frac{R_n( t -\frac{n}{2} + 1)}{R_n( t -\frac{n}{2} )}
= \frac{( t-\frac{n}{2} )^s\prod_{\theta \in \mathcal{Z}_B} ( t+\theta + \frac{n}{2} ) }{( t+1+\frac{n}{2} )^s\prod_{\theta \in \mathcal{Z}_B} ( t+\theta - \frac{n}{2} )}
= \frac{B_n(t)}{A_n(t)}
\end{equation}
with
\begin{align*}
A_n(z) &\coloneq \Big( z+1+\frac{n}{2}\Big)^{s+1}\Big(z-\frac{n}{2}\Big) \prod_{\theta \in \mathcal{Z}_B}\Big( z + \theta - \frac{n}{2} \Big), 
\\
B_n(z) &\coloneq \Big(z-\frac{n}{2}\Big)^{s+1}\Big( z+1+\frac{n}{2}\Big)\prod_{\theta \in \mathcal{Z}_B}\Big( z + \theta + \frac{n}{2} \Big). 
\end{align*}
In particular, the extension $a(x)$ by zero satisfies
\[
A_n(t)a(t+1)=B_n(t)a(t) \quad\text{for all}\; t\in \bbR,
\]
with polynomials $A_n, B_n$ of degree $d=s+|\mathcal{Z}_B|+2$. Since $\deg \widehat{R}_n = n|\mathcal{Z}_B| -(n+1)s$, we have
\begin{equation*} 
A_0\sum_{\substack{x \in \mathcal{X}_B \\ |x|>n/2}} \left| \widehat{w}_x \widehat{R}_n\left( x \right) \right| (1+|x|)^{(n+1)s-n|\mathcal{Z}_B|-2} < +\infty,
\end{equation*}
hence $\mu_n$ satisfies the growth condition \eqref{eq:growth-condition-GD} with $D=(n+1)s-n|\mathcal{Z}_B|-1$.
From \eqref{oddness_of_shifted_Rn} and the notation above, the number $\mathcal{L}_n(t^k)$ for $k<D$ defined in \eqref{def_L} is simply the $k$-th moment of the measure~$\mu_n$:
\begin{equation*} 
\mathcal{L}_n(t^k) = A_0\sum_{\substack{x \in \mathcal{X}_B \\ |x|>n/2}} \sgn(x)\widehat{w}_x \widehat{R}_n( x ) x^k.
\end{equation*}
We can now apply Proposition \ref{prop:nonvanishing-Hankel-det} to the signed measure $\mu_n$ to get

\begin{proposition}\label{prop_nonvanishing}
For any $n \in 2\bbZ_{>0}$ and $M \in \bbZ_{>0}$ with
\[
2 ( M + s +|\mathcal{Z}_B| +1 ) \le (n+1)s - n|\mathcal{Z}_B|,
\]
there exists $m \in \bbZ_{>0}$ satisfying
\[
M \le m \le M+ s +|\mathcal{Z}_B| +1 \quad\text{and}\quad \Delta_{n,m} \neq 0.
\]
\end{proposition}

\begin{proof}
By Proposition \ref{prop:nonvanishing-Hankel-det}, it suffices to show that there is no rational function $r(z)\in \bbC(z)$ satisfying the functional equation
\begin{equation}\label{eq:ArBrC}
A_n(z)r(z+1)-B_n(z)r(z)=C_n(z),
\end{equation}
with 
\[
C_n(z)\coloneq A_0\sum_{\substack{x \in \mathcal{X}_B \\ |x|>n/2}} \sgn(x)\widehat{w}_{x}\widehat{R}_n\left(x\right) \left(\frac{A_n(z)-A_n(x-1)}{z-(x-1)} - \frac{B_n(z)-B_n(x)}{z-x} \right).
\]
Suppose this is not true and there is a rational solution $r(z)$ of~\eqref{eq:ArBrC}. 
By Lemma \ref{lem:feq-CS-transform}, the Cauchy--Stieltjes transform
\begin{equation*} 
f_n(z) \coloneq A_0\sum_{\substack{x \in \mathcal{X}_B \\ |x|>n/2}} \frac{\sgn(x)\widehat{w}_x \widehat{R}_n( x )}{z-x}, \quad\text{where}\; z \in \bbC \setminus \big( \mathcal{X}_B \setminus [-n/2,n/2] \big),
\end{equation*}
of $\mu_n$ satisfies the same functional equation
\begin{equation}\label{AfBfC}
A_n(z) f_n(z+1) - B_n(z) f_n(z) = C_n(z).
\end{equation}
Subtracting \eqref{eq:ArBrC} from \eqref{AfBfC} and using
\eqref{5.8}, we find that the meromorphic function
\begin{equation*} 
g_n(z)=\frac{f_n(z)-r(z)}{\widehat{R}_n(z)}
\end{equation*}
satisfies
\begin{equation}\label{g_is_periodic}
g_n(z+1) = g_n(z).
\end{equation}
Now, take $b^* = \max\{ b \in \Psi_B : w_b \neq 0 \}$. 
Then, for any $T \in \bbZ_{>0}$, we have
\[
\widehat{w}_{T+1/b^*} = \widehat{w}_{-T+1/b^*} = \sum_{\substack{b \in \Psi_B \\ b^* \mid b}} w_b = w_{b^*} \neq 0.
\] 
Since $r(z)/\widehat{R}_n(z)$ is a rational function of~$z$, there is a choice of integer $T>n$ for which neither $T+1/b^*$ nor $-T+1/b^*$ is a pole of this rational function.
Then 
\[
\Res_{z=\pm T+1/b^*} g_n(z) 
= \frac1{\widehat{R}_n( \pm T+1/b^* )} \cdot \Res_{z=\pm T+1/b^*} f_n(z) = \pm A_0w_{b^*}.
\]
But \eqref{g_is_periodic} implies 
\[
\Res_{z=T+1/b^*} g_n(z) = \Res_{z=-T+1/b^*} g_n(z).
\]
which contradicts $A_0w_{b^*} \neq 0$. 
Therefore, equation \eqref{eq:ArBrC} does not have any rational solutions.
\end{proof}

\section{Proof of the main theorem}\label{sec6}

Fix any $\varepsilon \in (0,1)$ and take 
\begin{equation*} 
B = \bigg( \sqrt{\frac{4\zeta(6)}{\zeta(2)\zeta(3)}} - \frac{\varepsilon}{10} \bigg) \cdot \sqrt{\frac{s\log\log s}{\log s}}, \quad
C = \frac{s}{\log\log s}.	
\end{equation*}
By \eqref{est_Psi_B}, \eqref{est_Z_B}, \eqref{est_A_1}, and \eqref{est_A_2}, for any odd integer $s \ge s_0'(\varepsilon)$ we have
\[
s > 10B^2 > 100, \quad 10 < 2C < s-B^2-10, \quad |\mathcal{Z}_B| \le B^2-1,
\]
and
\begin{gather}
|\Psi_B| \ge \bigg( \sqrt{\frac{4\zeta(2)\zeta(3)}{\zeta(6)}} -\varepsilon \bigg) \cdot \sqrt{\frac{s \log\log s}{\log s}}, \label{6.2} \\
\log(A_1(B)A_2(B)) \le \Big( 1 -\frac{\varepsilon}{100} \Big) \cdot s\log\log s. \label{6.3}
\end{gather}
By applying Proposition \ref{prop_nonvanishing} with $M=\lceil Cn\rceil$, for any sufficiently large even integer $n$, there exists an integer $m=m(n)$ satisfying
\begin{equation}\label{property_of_m_n}
Cn \le m \le Cn + 2s \quad\text{and}\quad \Delta_{n,m} \neq 0;
\end{equation}
in particular, $m(n)/n \to C$ as $n \to \infty$.
By Proposition \ref{prop_ari} and Eq.~\eqref{property_of_m_n}, the number
\[
J_n = \frac{d_n^{sm}}{\prod_{i=0}^{m-1}( n!^{\lfloor i/n \rfloor} \cdot ( i - n\lfloor i/n \rfloor )! )^2} \cdot \Delta_{n,m}
\]
is a nonzero integer, so that $|J_n| \ge 1$. 
On the other hand, by Proposition \ref{prop_ana} and Eq.~\eqref{6.3},
\begin{align*}
\limsup_{n \to \infty} \frac{\log|J_n|}{n^2}
&\le C \log(A_1(B)A_2(B)) + C^2 + (2s+10B^2)C - \big( (s-B^2)C -2C^2 \big) \log \frac{C}{10B^2} \\
&= C \log(A_1(B)A_2(B)) - sC \log \frac{C}{10B^2} + o(s^2)
\\
&\le \Big(-\frac{\varepsilon}{100} + o(1)\Big) s^2
\quad\text{as}\; s \to \infty.
\end{align*}
This means that $\log|J_n|< 0$ for infinitely many $n$, whenever $s \ge s_0(\varepsilon)\ge s_0'(\varepsilon)$ is sufficiently large, a contradiction.
The source of the contradiction is our assumption \eqref{prop_I}, hence it cannot hold and our goal $N_s\ge |\Psi_B|$ is achieved.
Therefore, by \eqref{6.2} we deduce that
\[
N_s\ge |\Psi_B| \ge \bigg( \sqrt{\frac{4\zeta(2)\zeta(3)}{\zeta(6)}} -\varepsilon \bigg) \cdot \sqrt{\frac{s \log\log s}{\log s}}
\]
for any odd integer $s \ge s_0(\varepsilon)$.



\begin{thebibliography}{99}

\bibitem{ANL2021}
F. Adiceam, E. Nesharim and F. Lunnon,
\textit{On the $t$-adic Littlewood conjecture},
Duke Math. J. \textbf{170} (2021), no. 10, 2371--2419.

\bibitem{Ape1979}
R. Ap\'{e}ry,
\textit{Irrationalit\'{e} de $\zeta(2)$ et $\zeta(3)$},
Ast\'{e}risque \textbf{61} (1979), 11--13.




\bibitem{BR2001}
K. Ball and T. Rivoal,
\textit{Irrationalit\'{e} d'une infinit\'{e} de valeurs de la fonction z\^{e}ta aux entiers impairs},
Invent. Math. \textbf{146} (2001), no. 1, 193--207.

%
%
%
%
%
%
%
%


\bibitem{Bro2026}
F. Brown,
\textit{Mellin transforms, transfinite diameter and rational approximations of integrals},
Exp. Math. (2026), advance online publication, \href{https://doi.org/10.1080/10586458.2026.2695605}{doi:10.1080/10586458.2026.2695605}.


\bibitem{CDT2024+} F. Calegari, V. Dimitrov and Y. Tang, 
\textit{The linear independence of $1$, $\zeta(2)$, and $L(2,\chi_{-3})$}, 
Preprint (2024) \href{https://arxiv.org/abs/2408.15403}{{\tt arXiv:2408.15403}}.

%
%
%
%
%
%
%
%

\bibitem{DM2014}
D. Dominici and F. Marcell\'{a}n,
\textit{Discrete semiclassical orthogonal polynomials of class one},
Pacific J. Math. \textbf{268} (2014), no. 2, 389--411.

%
%


\bibitem{Fau2026+} A. Fauzan,
\textit{$\zeta(5)$ is irrational},
Preprint (2026) \href{https://doi.org/10.5281/zenodo.22826419}{zenodo:22826419}.	
%
%
%
%
%
%
%

\bibitem{Fis2026}
S. Fischler,
\textit{Linear independence of odd zeta values using Siegel's lemma},
J. London Math. Soc. (2) \textbf{113} (2026), no. 4, e70535, 42 pp.

%

\bibitem{FSZ2018}
S. Fischler, J. Sprang and W. Zudilin,
\textit{Many values of the Riemann zeta function at odd integers are irrational},
C. R. Math. Acad. Sci. Paris \textbf{356} (2018), no. 7, 707--711.

\bibitem{FSZ2019}
S. Fischler, J. Sprang and W. Zudilin,
\textit{Many odd zeta values are irrational},
Compos. Math. \textbf{155} (2019), no. 5, 938--952.

%

%
%
%
%



%
%



%
%
%
%


\bibitem{Lai2026a}
L. Lai,
\textit{A note on the number of irrational odd zeta values II},
J. Th\'{e}or. Nombres Bordeaux \textbf{38} (2026), no. 2, 469--487.

\bibitem{Lai2026b}
L. Lai,
\textit{Small improvements to the Ball--Rivoal theorem and its $p$-adic variant},
Front. Comb. Number Theory \textbf{4} (2026), 1--44.


\bibitem{LY2020}
L. Lai and P. Yu,
\textit{A note on the number of irrational odd zeta values},
Compos. Math. \textbf{156} (2020), no. 8, 1699--1717.


%
%
%
%
%
%
%
%


%
%
%
%

\bibitem{Riv2000}
T. Rivoal,
\textit{La fonction z\^{e}ta de Riemann prend une infinit\'{e} de valeurs irrationnelles aux entiers impairs},
C. R. Acad. Sci. Paris S\'{e}r. I Math. \textbf{331} (2000), no. 4, 267--270.

\bibitem{Riv2002}
T. Rivoal,
\textit{Irrationalit\'{e} d'au moins un des neuf nombres $\zeta(5),\zeta(7),\ldots,\zeta(21)$},
Acta Arith. \textbf{103} (2002), no. 2, 157--167.
%
%
%
%


%
%
%
%
%
%

\bibitem{Spr2018+}
J. Sprang,
\textit{Infinitely many odd zeta values are irrational. By elementary means},
Preprint (2018) \href{https://arxiv.org/abs/1802.09410}{arXiv:1802.09410}.



%
%
%
%
%
%
%
%
%
%
%
%
%
%

\bibitem{Zud2001b}
W. Zudilin,
\textit{One of the numbers $\zeta(5)$, $\zeta(7)$, $\zeta(9)$, $\zeta(11)$ is irrational},
Russian Math. Surveys \textbf{56} (2001), no. 4, 774--776.

%
%
%
%

\bibitem{Zud2004b}
W. Zudilin,
\textit{Arithmetic of linear forms involving odd zeta values},
J. Th\'{e}or. Nombres Bordeaux \textbf{16} (2004), no. 1, 251--291.

%
%
%

\bibitem{Zud2017}
W. Zudilin,
\textit{A determinantal approach to irrationality},
Constr. Approx. \textbf{45} (2017), no. 2, 301--310.

\bibitem{Zud2018}
W. Zudilin,
\textit{One of the odd zeta values from $\zeta(5)$ to $\zeta(25)$ is irrational. By elementary means},
SIGMA \textbf{14} (2018), 028, 8 pp.

\end{thebibliography}
\end{document}